\documentstyle[amssymb,amsfonts]{amsart}

\def\be#1{\begin{equation} \label{#1}}
\def\bi{\begin{itemize}}
\def\bs{\begin{split}}
\def\es{\end{split}}
\def\ba{\begin{align}}
\def\bas{\begin{align*}}
\def\ea{\end{align}}
\def\eas{\end{align*}}
\def\Im{{\hbox{Im}}}

\def\Re{{\hbox{Re}}}

\def\dist{{\hbox{\rm dist}}}

\def\D{{\mathcal D}}
\def\R{{\hbox{\bf R}}}

\def\T{{\hbox{\bf T}}}

\def\eps{\varepsilon}
\newenvironment{proof}{\noindent {\bf Proof} }{\endprf\par}
\def \endprf{\hfill  {\vrule height6pt width6pt depth0pt}\medskip}
\def\emph#1{{\it #1}}
\def\textbf#1{{\bf #1}}
\def\divider#1{$\bullet\quad${\bf #1}}

\theoremstyle{plain}
  \newtheorem{theorem}[subsection]{Theorem}

  \newtheorem{lemma}[subsection]{Lemma}

\theoremstyle{remark}

\theoremstyle{definition}

\numberwithin{equation}{section}

\title[Stability bounds below the energy norm]
{Polynomial upper bounds for the orbital\\ instability of the 1D
cubic NLS\\ below the energy norm}

\author{J.~Colliander}
\thanks{J.E.C. is supported in part by N.S.F.
 grant DMS 0100595 and N.S.E.R.C. grant RGPIN 250233-03.}
\address{University of Toronto}

\author{M.~Keel}
\thanks{M.K. is supported in part by N.S.F. Grant DMS
                         9801558}
\address{University of Minnesota}

\author{G.~Staffilani}
\thanks{G.S. is supported in part by N.S.F.
Grant DMS 0100375 and by a grant from  the Sloan Foundation.}
\address{Massachusetts Institute of Technology}

\author{H.~Takaoka}
\address{Hokkaido University}
\thanks{H.T. is supported in part by J.S.P.S. Grant No. 13740087.}

\author{T.~Tao}
\thanks{T.T. is a Clay Prize Fellow and is supported in part by a grant
from the Packard Foundation.}
\address{University of California, Los Angeles}


\begin{document}

\maketitle

\begin{abstract}
We study the long-time behaviour of the focusing cubic NLS on $\R$ in the Sobolev norms $H^s$ for $0 < s < 1$.  We obtain polynomial growth-type upper bounds on the $H^s$ norms, and also limit any orbital $H^s$ instability of the ground state to polynomial growth at worst; this is a partial analogue of the $H^1$ orbital stability result of Weinstein \cite{wein}, \cite{wein:modulate}.  In the sequel to this paper we generalize this result to other nonlinear Schr\"odinger equations.  Our arguments are based on the ``$I$-method'' from earlier papers \cite{ckstt:1}-\cite{ckstt:8} which pushes down from the energy norm, as well as an ``upside-down $I$-method'' which pushes up from the $L^2$ norm.
\end{abstract}

\maketitle

\section{Introduction}\label{introduction-sec}

We consider the long-time behaviour of solutions to the Cauchy problem for the one-dimensional focusing cubic nonlinear Schr\"odinger equation
\be{nls}
iu_t + u_{xx} = -F(u); \quad u(x,0) = u_0(x)
\end{equation}
where $u(x,t)$ is a complex-valued function on $\R \times \R$, $F(u)$ is the focusing cubic nonlinearity $F(u) := u \overline{u} u$, and $u_0(x)$ lies in the Sobolev space $H^s(\R)$ for some $s \in \R$.

It is known (\cite{tsutsumi}, \cite{cwI}) that the Cauchy problem \eqref{nls} is globally well-posed in $H^s$ for all\footnote{For $s < 0$ one does not even have local well-posedness, at least if one demands uniform continuity of the solution map; see \cite{kpv:counter}.} $s \geq 0$.  Furthermore, due to the many conservation laws of \eqref{nls}, we know that if $s$ is an integer and the initial data is in $H^s$, then the $H^s$ norm stays bounded for all time.
For instance, for $s=1$ one can obtain uniform $H^1$ bounds by exploiting the conservation of the Hamiltonian
$$ H(u) := \int  \frac{1}{2} |u_x|^2 - \frac{1}{4} |u|^4\ dx$$
and the $L^2$ norm, combined with the Gagliardo-Nirenberg inequality.

However when $s$ is not an integer, the standard iteration argument only gives bounds on the $H^s$ norms which grow exponentially in time.  For $s>1$ there are
polynomial growth bounds\footnote{It is an open question whether one has some sort of scattering in $H^s$ for this equation, which would of course imply that the $H^s$ norm remains bounded.  Some recent progress in this direction is in \cite{nak:scatter}.} in \cite{staff:growth}, \cite{borg:growth}, \cite{cdks}.

The first result of this paper is to extend these techniques to $s$ between 0 and 1.

\begin{theorem}\label{poly-growth}  If $0 < s < 1$ and $u_0 \in H^s$, then we have
$$ \| u(t) \|_{H^s} \leq C(\|u_0\|_{H^s}) (1 + |t|)^{2s+}.$$
\end{theorem}

The proof of this theorem proceeds by an ``upside-down'' version of the
``$I$-method'' (\cite{ckstt:1}, \cite{ckstt:3}, \cite{ckstt:4}, \cite{ckstt:5},
\cite{ckstt:6}, \cite{keel:mkg}; see also \cite{borg:book}, \cite{ckstt:2},
\cite{keel:wavemap}), in which one applies a differentiation operator $\D = \D_N$
to the solution instead of a smoothing operator $I = I_N$, and proves that the
quantity $\| \D_N u(t) \|_2^2$
is almost conserved in time. One should compare the result in Theorem
\ref{poly-growth} with the ones established in
\cite{borg:growth}, \cite{cdks} and \cite{staff:growth} for the time
asymptotic  of the $H^{s}$ norm  for smooth ($s>1$) KdV and
Schr\"odinger type solutions. The method used
by Bourgain in \cite{borg:growth} is also based on an improved local
estimate, but the improvement is not obtained by replacing the
$H^{s}$ norm with a better integral (which is our argument here), but
instead by using well-posedness results below the energy norm $H^{1}$.
In \cite{cdks} and \cite{staff:growth} the improvement of the local
estimate is obtained by using sharp bilinear estimates in negative
Sobolev spaces.

We remark that the same polynomial growth result in Theorem
\ref{poly-growth} also holds for the defocusing
cubic NLS (in which $-|u|^2 u$ is replaced by $+|u|^2 u$) and is slightly easier
to obtain.
It is likely that one can use the correction term techniques in \cite{ckstt:2},
\cite{ckstt:6} to improve the exponent $2s+$ substantially\footnote{In the
particular case of the cubic 1D NLS equation \eqref{nls}, one may also be
able to exploit the complete integrability of the equation to obtain bounds
on the $H^s$ norms which are uniform in time, and perhaps also to obtain
global stability bounds for solitons and multisolitons as well.  On the
other hand, the methods here do not exploit complete integrability and
are applicable to a wide range of Hamiltonian evolution equations.},
perhaps all the way down to $0+$.  Certainly one expects to obtain
an exponent which goes to 0 as $s \to 1^-$ by exploiting conservation
of the Hamiltonian.

One can view Theorem \ref{poly-growth} as a bound on the possible $H^s$ instability of the 0 solution $u \equiv 0$ to \eqref{nls}.  It is not strong enough to say that small $H^s$ perturbations to this solution at time zero remain small perturbations for all later time, but it limits the growth of the perturbation to polynomial growth at worst.

The next result of this paper concerns the $H^s$ orbital stability of ground states for \eqref{nls}.  For simplicity we shall only consider the ground states at energy 1 (the other energies can then be recovered by a scaling argument\footnote{More specifically, for every energy $E > 0$, there is a unique positive even Schwartz function $Q_E$ obeying $(Q_E)_{xx} + |Q_E|^2 Q_E = EQ_E$, but these ground states are linked by the scaling $Q_E(x) = E^{1/2} Q(E^{1/2} x)$.  Because the equation \eqref{nls} is $L^2$-subcritical, all of these ground states have different $L^2$ mass.  Since the $L^2$ mass is an invariant of the NLS flow, we can thus restrict our attention to a sphere in $L^2$, in which case only one energy $E$ is relevant.  One can then use the scale invariance of \eqref{nls} to set $E=1$.}). It is known \cite{coff} that there exists a unique even positive Schwartz function $Q(x)$ on $\R$ which solves the equation\footnote{Indeed, we have the explicit formula $Q(x) := 2^{-1/2} / \cosh(x)$, although we will not use this formula in this paper.}
\be{Q-def}
Q_{xx} + |Q|^2 Q = Q.
\end{equation}
The Cauchy problem \eqref{nls} with initial data $u_0 = Q$ then has an explicit solution $u(t) = e^{it} Q$.  More generally, for any $x_0 \in \R$ and $e^{i\theta} \in S^1$, the Cauchy problem with initial data $u_0(x) = e^{i\theta} Q(x-x_0)$ has explicit solution $e^{i(\theta+t)} Q(x-x_0)$.  If we thus define the two-dimensional \emph{ground state cylinder}\footnote{Note that the ground state cylinder is the orbit of $Q$ under the phase and translation invariances of NLS.  We do not utilize the scaling invariance because, as mentioned earlier, this changes the $L^2$ norm of $Q$.  Also we do not utilize Gallilean invariance because this does not preserve the Hamiltonian.}  $\Sigma \subset H^1(\R)$ by
$$ \Sigma := \{ e^{i\theta} Q(\cdot - x_0): x_0 \in \R, e^{i\theta} \in S^1 \}$$
we see that the nonlinear flow \eqref{nls} preserves $\Sigma$.  Also note that each element of $\Sigma$ obeys \eqref{Q-def} (though of course most ground states are not even or positive).

In \cite{wein} (see also \cite{wein:modulate}) Weinstein showed that the
ground state cylinder $\Sigma$ was $H^1$-stable.  More precisely, he showed
an estimate of the form
\be{stab}
\dist_{H^1}(u(t), \Sigma) \sim \dist_{H^1}(u_0, \Sigma)
\end{equation}
(when $\dist_{H^1}(u_0, \Sigma)$ is small),
for all $H^1$ solutions $u(t)$ to \eqref{nls} and all times $t \in \R$.  In other words, solutions which start close to a ground state in $H^1$ at time $t=0$, will stay close to a ground state for all time (though the nearby ground state may itself vary in time\footnote{For instance, consider the solution $u(x,t) = e^{i(\eps x-\eps^2 t)} e^{it} Q(x-2\eps t)$ for some small $\eps$; this is a Gallilean transformation of the ground state solution $u(t) = e^{it} Q$, and is close to this solution at time zero.  However at later times, the solution slowly drifts away from the original ground state solution, although it remains close to the ground state cylinder $\Sigma$.}).

To prove \eqref{stab}, Weinstein employed the \emph{Lyapunov functional}\footnote{This is the functional for energy $E=1$.  For other energies it is given by $L(u) = 2H(u) + \int E |u|^2$.}
\be{l-def}
L(u) := 2H(u) + \int |u|^2
= \int |u_x|^2 + |u|^2 - \frac{1}{2} |u|^4\ dx,
\end{equation}
which is well-defined for all $u \in H^1$.  Since this quantity is a combination of the Hamiltonian and the $L^2$ norm, it is clearly an invariant of the flow \eqref{nls}.  More explicitly, for sufficiently smooth functions $u(x,t)$ we have the formula
\be{L-deriv}
\partial_t L(u)
= 2 \langle u_t, -u_{xx} + u - F(u)\rangle
\end{equation}
which clearly vanishes if $u$ solves \eqref{nls}.  Here and in the sequel we use $\langle, \rangle$ to denote the \emph{real} inner product
$$ \langle u, v \rangle := \Re \int u \overline{v}\ dx.$$

>From \eqref{L-deriv} and \eqref{Q-def} we see that the ground states in $\Sigma$ are critical points of $L$.  In fact they are minimizers of $L$; more precisely, we have the fundamental \emph{coercivity estimate}
\be{weinstein}
L(u) - L(Q) \sim \dist_{H^1}(u, \Sigma)^2 \hbox{ whenever } u \in H^1 \hbox{ and } \dist_{H^1}(u, \Sigma) \ll 1;
\end{equation}
see \cite{wein}.  The stability estimate \eqref{stab} then follows easily from \eqref{weinstein} and the conservation of $L$.

Note that the functional $L$ is invariant under phase rotation $u \to e^{i\theta} u$ and translation $u \mapsto u(\cdot - x_0)$.  Thus one cannot expect a coercivity bound like \eqref{weinstein} in these directions.  Of course, this is consistent with \eqref{weinstein} since the ground state cylinder $\Sigma$ is itself invariant under these symmetries; the point of \eqref{weinstein} is that there are no other directions (in the tangent space of $H^1$ at $Q$) for which $L$ can be invariant or concave.

The second main result of this paper is to partially extend the $H^1$ orbital
stability result to an $H^s$ orbital stability-type result for $0 \leq s < 1$.  Unfortunately, as in Theorem \ref{poly-growth}, our estimate loses a polynomial factor in $t$, so we cannot exclude the possibility of polynomial orbital instability in $H^s$:

\begin{theorem}\label{main-2}  Let $0 \leq s < 1$, and suppose $\dist_{H^s}(u_0, \Sigma) \ll 1$.  Then we have
\be{stab-s}
\dist_{H^s}(u(t), \Sigma) \lesssim t^{1-s+} \dist_{H^s}(u_0, \Sigma)
\end{equation}
whenever
 $ 1 \leq t \ll \dist_{H^s}(u_0, \Sigma)^{-1/(1-s+)}.$
In particular, $u(t)$ stays in a bounded subset of $H^s$ for all
times $|t| \ll \dist_{H^s}(u_0, \Sigma)^{-1/(1-s+)}$.
\end{theorem}

Note that a naive application of the local well-posedness theory would lose a factor of the form $\exp(C t)$ in the estimate \eqref{stab-s}, so that one could only assure $u(t)$ stays in a bounded subset of $H^s$ for times $|t| \lesssim \log (1/\dist_{H^s}(u_0, \Sigma))$.  Note that for times $t$ close
to $\dist_{H^s}(u_0, \Sigma)^{1/(1-s+)}$, the right-hand side of \eqref{stab-s} is close to 1, so that we are no longer keeping $u(t)$ close to $\Sigma$.  At this point one can use Theorem \ref{poly-growth} to control the further development of $u(t)$.

The proof of \eqref{stab-s} proceeds via the smoothing operator $I = I_N$ mentioned earlier.  Basically, the idea is to show that the modified Lyapunov functional $L(Iu)$ is almost conserved, and then combine this with \eqref{weinstein} to obtain \eqref{stab-s}.  It turns out that a naive implementation of this approach loses an epsilon power of $\dist_{H^s}(u_0, \Sigma)$ because the operator $I$ does not quite preserve the ground state cylinder, but this can be rectified by the standard technique of choosing an approximating ground state to obey specially chosen orthogonality conditions.  For expository reasons we have chosen to give the naive versions of the argument first, and only give the full argument at the end of the paper.

The factor $t^{1-s+}$ in \eqref{stab-s} can probably be reduced, however the factor $\dist_{H^s}(u_0, \Sigma)$ on the right-hand side is necessary (as can be seen even when $t \sim 1$).

In all of these arguments it is crucial that the $L^2$ and $H^1$ norms are both subcritical (in the sense that they scale as a negative power of length using the natural scaling of \eqref{nls}).  It is because of this that we cannot extend these results to any other NLS with an algebraic nonlinearity\footnote{More specifically, the NLS equation is only $L^2$ subcritical when $p < 1 + \frac{4}{n}$, while an algebraic nonlinearity only occurs when $p$ is an odd integer.  Since $p>1$, the only subcritical algebraic equation occurs when $n=1$, $p=3$.}.  In the sequel to this paper we shall obtain some partial results of the above type in the case when $p$ is not an odd integer; the main new difficulty is to commute the $I$ operator with non-algebraic nonlinearities $F(u)$.

\section{Notation}\label{notation-sec}

We use $A \lesssim B$ to denote $A \leq CB$, where $C$ is a constant depending on $s$ which may vary from line to line.  We use $a+$, $a-$ to denote quantities of the form $a+\eps$, $a-\eps$, where $\eps$ is arbitrarily small.  We use $\langle \xi \rangle$ to denote $1 + |\xi|$.

We define the spatial Fourier transform by
$$ \hat f(\xi) := \int_\R e^{-i x \xi} f(x)\ dx$$
and the spacetime Fourier transform by
$$ \tilde u(\xi,\tau) := \int_\R \int_\R e^{-i (x \xi+t\tau)} u(x,t)\ dx dt.$$
Following \cite{borg:xsb}, we define the $X^{s,b}$ spaces by the norm
$$ \| u \|_{X^{s,b}} := \| u \|_{s,b} := \| \langle \xi \rangle^s \langle \tau - \xi^2 \rangle^b \tilde u(\xi,\tau) \|_{L^2_\tau L^2_\xi}.$$
For any time interval $I$, we define the restricted $X^{s,b}$ spaces by
$$ X^{s,b}_I := \{ u|_{\R \times I}: u \in X^{s,b} \}$$
with the usual norm
$$ \| v \|_{X^{s,b}_I} := \inf \{ \| u \|_{s,b}: u|_{\R \times I} = v \}.$$

We shall need the Strichartz estimate
\be{strichartz-6}
\| u \|_{L^6_{x,t}} \lesssim \| u \|_{0,1/2+};
\end{equation}
this can be obtained by writing an $X^{0,1/2+}$ function as an average of modulated free Schr\"odinger waves (as in \cite{ginibre:survey}) and then using the $L^6_{x,t}$ Strichartz estimate for free solutions (see e.g. \cite{tao:keel} and the references therein).

Interpolating \eqref{strichartz-6} with the trivial bound $\| u \|_{L^2_{x,t}} \leq \| u\|_{0,0}$ we obtain
\be{strichartz-4}
\| u \|_{L^4_{x,t}} \lesssim \| u \|_{0,3/8+}.
\end{equation}
We also record the variant estimate
\be{strichartz-84}
\| u \|_{L^8_t L^4_x} \lesssim \| u \|_{0,1/2+};
\end{equation}
this can be obtained by interpolating \eqref{strichartz-6} with the energy estimate $\| u \|_{L^\infty_t L^2_x} \lesssim \| u \|_{0,1/2+}$.

>From \cite{borg:book} (see also \cite{ckstt:5}) we recall the improved bilinear Strichartz estimate (in one spatial dimension)
\be{improved-strichartz}
\| u_1 u_2 \|_{L^2_{t,x}} \lesssim N^{-1/2} \| u_1\|_{0,1/2+} \|u_2\|_{0,1/2+}
\end{equation}
whenever $u_1$ has Fourier support in the region $|\xi| \sim N$, and $u_2$ has Fourier support in the region $|\xi| \ll N$.

Let $n \geq 2$, and let $m(\xi_1, \ldots, \xi_n)$ be a function supported on the hyperplane $\{ \xi_1 + \ldots + \xi_n = 0 \}$.  We use the notation
$$ \Lambda_n(m(\xi_1, \ldots, \xi_n); f_1, \ldots, f_n)$$
to denote the multilinear form
$$ \Lambda_n(m(\xi_1, \ldots, \xi_n); f_1, \ldots, f_n) := \int_{\xi_1 + \ldots + \xi_n = 0} m(\xi_1, \ldots, \xi_n) \hat f_1(\xi_1) \ldots \hat f_n(\xi_n).$$

\section{Proof of Theorem \ref{poly-growth}}\label{growth-sec}

We shall divide the proof of this theorem into several broad steps.  These steps will also appear in the proof of Theorem \ref{main-2}.

\medskip

\divider{Step 0.  Preliminaries; introduction of the modified energy.}

Fix $0 < s < 1$, and fix the global $H^s$ solution $u$.  Henceforth all implicit constants may depend on $s$ and the quantity $\| u_0 \|_{H^s}$.  Our task is to show that
\be{t-grow}
\| u(T)\|_{H^s} \lesssim T^{2s+}
\end{equation}
for all $T \gg 1$.  By the $H^s$ global well-posedness and regularity theory (see e.g. \cite{cwI}) and the usual limiting argument it suffices to do this for smooth, rapidly decreasing $u$.

We now apply an ``upside-down'' version of the $I$-method; whereas the strategy of the $I$ method is to mollify the solution at high frequencies to make it smoother (e.g. in the energy class $H^1$), here we amplify the solution at high frequencies instead to make it rougher (specifically, we place it in $L^2$).  In contrast, the proof of Theorem \ref{main-2} in later sections will proceed via the ordinary ``$I$-method''.

Fix $T$, and let $N \gg 1$ be a large quantity depending on $T$ to be chosen later.
Let $\theta(\xi)$ be a smooth even real-valued symbol such that
$\theta(\xi) = 1$ for $|\xi| \leq N$ and $\theta(\xi) =
|\xi|^s/N^s$ for $|\xi| > 2N$, and let $\D$ be the Fourier multiplier
\be{d-def}
\widehat{\D f}(\xi) := \theta(\xi) \hat f(\xi).
\end{equation}
Thus $\D$ is the identity for low frequencies $|\xi| \leq N$, and becomes a differentiation operator of order $s$ for high frequencies $|\xi| \gtrsim N$.

Define the \emph{modified energy} $E_N(t)$ at time $t$ by
\be{en-def}
E_N(t) := \| \D u(t)\|_2^2.
\end{equation}
>From Plancherel we have the upper bound
$$ E_N(0) \lesssim \| u_0 \|_{H^s}^2 \lesssim 1.$$

The heart of the argument shall lie in the following almost conservation law for $E_N$.

\begin{lemma}\label{a-c}  If $t_0 \in \R$ is such that $E_N(t_0) \leq C$ for some bounded constant $C = O(1)$, then we have
\be{ftoc}
E_N(t_0+\delta) = E_N(t_0) + O(N^{-1/2+})
\end{equation}
where $\delta > 0$ is an absolute constant depending only on $s$ and $C$.
\end{lemma}

The error bound of $O(N^{-1/2+})$ might not be sharp; any improvement in this estimate will lead to a better polynomial growth bound than \eqref{t-grow}.  It may be that one can improve this result by adding suitably chosen ``correction terms'' to $E_N(t)$, in the spirit of \cite{ckstt:2}, \cite{ckstt:6}.

\medskip

\divider{Step 1.  Deduction of \eqref{t-grow} from Lemma \ref{a-c}.}

If we assume Lemma \ref{a-c}, then we may iterate it to obtain
$E_N(t) \lesssim 1$ for all $0 \leq t \ll N^{1/2-}$.  In particular, if we assume $T \ll N^{1/2-}$, we have
from Plancherel and \eqref{en-def} that
$$ \| u(T) \|_{H^s} \lesssim N^s E_N(T)^{1/2} \lesssim N^s.$$
Optimizing $N$ in terms of $T$ we obtain \eqref{t-grow} as desired.

It remains to prove Lemma \ref{a-c}.  This is done in several stages.

\medskip

\divider{Step 2.  Control $u$ at time $t_0$.}

By hypothesis we have
\be{init}
\| \D u(t_0) \|_2 \lesssim 1.
\end{equation}
As one can see from \eqref{d-def}, this is essentially an $H^s$-type bound on $u$ (up to powers of $N$).  However, in order to obtain good polynomial growth bounds it is important that we use the norm $\| \D u \|_2$ throughout rather than $\| u \|_{H^s}$.

\medskip

\divider{Step 3.  Control $u$ on the time interval $[t_0-\delta, t_0 + \delta]$.}

We now use \eqref{init} to claim

\be{spacetime}
\| \D u \|_{X^{0,1/2+\eps}_{[t_0-\delta,t_0+\delta]}} \lesssim 1
\end{equation}
for any $0 < \eps \ll 1$,
if $0 < \delta \ll 1$ is a sufficiently small constant (depending on $s$ and the bound in \eqref{init}, but not on $N$).  This will be achieved by techniques similar to those used to obtain local well-posedness using the $X^{s,b}$ spaces (as in e.g. \cite{borg:xsb}).

To begin with, by the standard energy estimate for $X^{s,b}$ spaces (see e.g. \cite{borg:xsb}, or \cite{cst:kdv}, page 771; note that the multiplier $\D$ commutes with the Schr\"odinger operator $i \partial_t + \partial_{xx}$ and is thus harmless) we have
$$ \| \D u \|_{X^{0,1/2+\eps}_{[t_0-\delta,t_0+\delta]}} \lesssim
\| \D u(t_0) \|_{L^2} + \delta^{0+} \| \D(iu_t + u_{xx}) \|_{X^{0,-1/2+2\eps}_{[t_0-\delta,t_0+\delta]}}.$$
Since we are allowed to choose $\delta$ to be sufficiently small, it suffices from \eqref{init}, \eqref{nls} and standard continuity (or iteration) arguments to show the trilinear estimate
$$ \| \D(u_1 \overline u_2 u_3) \|_{0,-1/2+2\eps} \lesssim \| \D u_1 \|_{0, 1/2+} \| \D u_2 \|_{0,1/2+} \| \D u_3 \|_{0,1/2+}.$$
We may assume that the $u_i$ have non-negative Fourier transforms.  Since $w(\xi + \eta) \lesssim w(\xi) + w(\eta)$ we can obtain a fractional Leibnitz rule, and reduce to showing
\be{threesome}
\begin{split}
 \| \D(u_1) \overline u_2 u_3 \|_{0,-1/2+2\eps} &\lesssim \| \D u_1 \|_{0, 1/2+} \| \D u_2 \|_{0,1/2+} \| \D u_3 \|_{0,1/2+}\\
\| u_1 \D(\overline u_2) u_3 \|_{0,-1/2+2\eps} &\lesssim \| \D u_1 \|_{0, 1/2+} \| \D u_2 \|_{0,1/2+} \| \D u_3 \|_{0,1/2+}\\
\| u_1 \overline u_2 \D(u_3) \|_{0,-1/2+2\eps} &\lesssim \| \D u_1 \|_{0, 1/2+} \| \D u_2 \|_{0,1/2+} \| \D u_3 \|_{0,1/2+}.
\end{split}
\end{equation}
>From the dual of \eqref{strichartz-4} we see that $\| f \|_{0,-1/2+2\eps} \lesssim \|f \|_{L^{4/3}_{x,t}}$.  The claims \eqref{threesome} then follow from H\"older's inequality and several applications of \eqref{strichartz-4}.  This gives \eqref{spacetime}.

 \medskip

\divider{Step 4.  Control the increment of the modified energy.}

To prove \eqref{ftoc}, we now apply the fundamental theorem of Calculus to write the increment of the modified energy \eqref{en-def} as
\be{the-ftoc}
E_N(t_0+\delta) - E_N(t_0) =  \int_{t_0}^{t_0+\delta} \partial_t E_N(t)\ dt.
\end{equation}
A routine integration by parts\footnote{This is of course related to the usual energy method computations to control the growth of higher order energies; see also \cite{borg:growth}, \cite{staff:growth}, \cite{cdks}.  Observe that if $\D$ was the identity then this derivative would vanish, thus giving the standard proof of $L^2$ norm conservation.} shows that
\bas
\partial_t E_N(t) &= 2 \langle \D u_t, \D u \rangle \\
&= 2 \langle i\D u_{xx} + i\D F(u),  \D u\rangle\\
&= 2 \langle i u \overline u u, \D^2 u \rangle \\
&= -2 \Im \Lambda_4(\theta(\xi_4)^2; u(t), \overline{u(t)}, u(t), \overline{u(t)})\\
&= \frac{1}{2} \Im \Lambda_4(\theta(\xi_1)^2 - \theta(\xi_2)^2 + \theta(\xi_3)^2
- \theta(\xi_4)^2; u(t), \overline{u(t)}, u(t),
\overline{u(t)}),
\end{align*}
where in the last step we exploited the symmetry
$$
\Lambda_4(\theta(\xi_j)^2; u(t), \overline{u(t)}, u(t), \overline{u(t)})
= \overline{\Lambda_4(\theta(\xi_k)^2; u(t), \overline{u(t)}, u(t),
\overline{u(t)})}$$
for $j=1,3$ and $k = 2,4$.

>From the above computations, it suffices to show the estimate
$$ |\int \chi_{[t_0,t_0+\delta]}(t)
\Lambda_4(M_4; u_1(t), \overline{u_2(t)}, u_3(t), \overline{u_4(t)})\ dt|
\lesssim N^{-1/2+} \prod_{i=1}^4 \| \D u_i \|_{0,1/2+}
$$
for all functions $u_1, u_2, u_3, u_4 \in X^{0,1/2+}$, where
$M_4 = M_4(\xi_1,\xi_2,\xi_3,\xi_4)$ denotes the symbol
$$ M_4(\xi_1,\xi_2,\xi_3,\xi_4) := \theta(\xi_1)^2 - \theta(\xi_2)^2 +
\theta(\xi_3)^2 - \theta(\xi_4)^2.$$

We shall argue similarly to \cite{ckstt:5}.  We apply Littlewood-Paley  decompositions, and assume that $u_i$ is supported on the region $\langle \xi_i \rangle \sim N_i$ for some dyadic $N_i \geq 1$; of course, we will eventually have to sum in $N_i$ to recover the general case.  Define $\{soprano, alto, tenor, baritone\} = \{1,2,3,4\}$ by requiring
$$ N_{soprano} \geq N_{alto} \geq N_{tenor} \geq N_{baritone}.$$
We may assume that $N_{soprano} \sim N_{alto}$ since $\Lambda_4$ is integrating over the region $\xi_1 + \xi_2 + \xi_3 + \xi_4 = 0$.  We may also assume that $N_{soprano} \gtrsim N$ since the symbol
$ \theta(\xi_1)^2 - \theta(\xi_2)^2 + \theta(\xi_3)^2 - \theta(\xi_4)^2$
vanishes otherwise.

We divide into two cases.
 \medskip

\divider{Case (a): $N_{soprano} \gg N_{baritone}$.}

In this case $M_4 = O(\theta(N_{soprano})^2)$, so we can essentially
\footnote{To be more precise, one should replace $u_i$ by $Mu_i$, the
Hardy-Littlewood maximal  function of $u_i$.  This is because the symbol
$M_4$ is equal to $\theta(N_{soprano})^2$ times a smooth multiplier which obeys
good kernel estimates.  We omit the details.  An alternate approach is to
reduce to the case when all the $ u_i$ have non-negative spacetime Fourier
transform (which is legitimate since the $X^{s,b}$ norms do not care about the phase of $\tilde u_i$), however the time cutoff $\chi_{[t_0, t_0+\delta]}(t)$ then presents an annoying technical difficulty (since its Fourier transform is not non-negative).} bound this term by
$$ \theta(N_{soprano})^2 \int\int_{t_0}^{t_0+\delta} |u_{soprano}| |u_{alto}|
|u_{tenor}| |u_{baritone}|\ dx dt.$$
We use H\"older's inequality to take $u_{alto}$ and $u_{tenor}$ in
$L^4_{x,t}$, and $u_{soprano} u_{baritone}$ in $L^2_{x,t}$.
Using two applications of \eqref{strichartz-4} and one application of
\eqref{improved-strichartz} we can bound this by
$$ \theta(N_{soprano})^2 N_{soprano}^{-1/2}
\prod_{i=1}^4 \| u_i \|_{0,1/2+}
\lesssim N_{soprano}^{-1/2} \prod_{i=1}^4 \| \D u_i \|_{0,1/2+}.$$
The claim then follows by summing in the $N_i$.

\medskip

\divider{Case (b): $N_{soprano} \sim N_{baritone}$.}

In this case all the frequencies are comparable.  The key observation is that
for $|\xi_i| \sim N_{soprano}$ and $\xi_1 + \xi_2 + \xi_3 + \xi_4 = 0$, we have
\be{filter}
\bigl|\frac{ M_4 }
{ \xi_1^2 - \xi_2^2 + \xi_3^2 - \xi_4^2 }\bigr| \lesssim
\frac{\theta(N_{soprano})^2}{N_{soprano}^2}.
\end{equation}
To see this, we write the denominator as
\bas
\xi_1^2 - \xi_2^2 + \xi_3^2 - \xi_4^2 &= (\xi_1-\xi_2)(\xi_1+\xi_2) + (\xi_3-\xi_4)(\xi_3+\xi_4)\\
&= (\xi_1-\xi_2-\xi_3+\xi_4)(\xi_1 +\xi_2)\\
&= 2(\xi_1+\xi_4)(\xi_1+\xi_2)
\end{align*}
and the numerator as
\bas
M_4 = &\theta(\xi_1)^2 - \theta(\xi_2)^2 + \theta(\xi_3)^2 - \theta(\xi_4)^2 \\
= &\theta(\xi_3 + (\xi_1+\xi_2) + (\xi_1+\xi_4))^2
- \theta(\xi_3 + (\xi_1+\xi_4))^2 \\
&+
\theta(\xi_3)^2 - \theta(\xi_3 + (\xi_1 + \xi_2))^2;
\end{align*}
the claim then follows from the double mean value theorem and the estimate
$$ |\frac{d^2}{d\xi^2} (\theta(\xi)^2)| \lesssim \frac{\theta(N_{soprano})^2}
{N_{soprano}^2}$$
for all $\xi = O(N_{soprano})$.

To use \eqref{filter} we must use the spacetime Fourier transform, which requires us to first modify the cutoff $\chi_{[t_0,t_0+\delta]}(t)$.  Write $\chi_{[t_0,t_0+\delta]}(t) = a(t) + b(t)$, where $a(t)$ is $\chi_{[t_0,t_0+\delta]}(t)$ convolved with a smooth approximation to the identity of width $N_{soprano}^{-100}$, and $b(t) = \chi_{[t_0,t_0+\delta]}(t) - a(t)$.

Consider the contribution of $b(t)$.  This term is essentially
bounded by
$$ \theta(N_{soprano})^2 \int\int |b(t)| |u_1(t)| |u_2(t)| |u_3(t)| |u_4(t)|\ dx dt;$$
by H\"older's inequality and \eqref{strichartz-84} we may bound this by
$$ N_{soprano}^{-10} \prod_{i=1}^4 \| u_i \|_{0,1/2+}$$
which easily sums to be acceptable.

Now consider the contribution of $a(t)$.  In light of the estimate
$$ \| a(t) u_1 \|_{0,1/2+} \lesssim N_{soprano}^{0+} \| u_1 \|_{0,1/2+}$$
(see \cite{ckstt:5}) it suffices to show
$$ \int
\Lambda_4(M_4; u_1(t), \overline{u_2(t)}, u_3(t), \overline{u_4(t)})\ dt
\lesssim N_{soprano}^{-1} \prod_{i=1}^4 \| \D u_i \|_{0,1/2+}.
$$
We may assume that the spacetime Fourier transforms of $u_i$ are non-negative.  By Plancherel and \eqref{filter} it suffices to show that
\bas
\int_{\xi_1 + \xi_2 + \xi_3 + \xi_4 = 0} \int_{\tau_1 + \tau_2 + \tau_3 + \tau_4 = 0} |\xi_1^2 - \xi_2^2 + \xi_3^2 - \xi_4^2| &
\tilde u_1(\xi_1,\tau_1) \tilde{\overline{u_2}}(\xi_2,\tau_2)
\tilde u_3(\xi_3,\tau_3) \tilde{\overline{u_4}}(\xi_4,\tau_4) \\
&\lesssim N_{soprano} \prod_{i=1}^4 \| u_i \|_{0,1/2+}.
\end{align*}
>From the triangle inequality we have
$$ |\xi_1^2 - \xi_2^2 + \xi_3^2 - \xi_4^2| \lesssim |\tau_1 - \xi_1^2| + |\tau_2 + \xi_2^2| + |\tau_3 - \xi_3^2| + |\tau_4 + \xi_4^2|;$$
interpolating this with the trivial bound of $N_{soprano}^2$ we obtain
$$ |\xi_1^2 - \xi_2^2 + \xi_3^2 - \xi_4^2| \lesssim N_{soprano} (|\tau_1 - \xi_1^2|^{1/2} + |\tau_2 + \xi_2^2|^{1/2} + |\tau_3 - \xi_3^2|^{1/2} + |\tau_4 + \xi_4^2|^{1/2}).$$
By symmetry, it thus suffices to show
\bas |\int_{\xi_1 + \xi_2 + \xi_3 + \xi_4 = 0} \int_{\tau_1 + \tau_2 + \tau_3 + \tau_4 = 0} |\tau_1 - \xi_1^2|^{1/2} &
\tilde u_1(\xi_1,\tau_1) \tilde{\overline{u_2}}(\xi_2,\tau_2)
\tilde u_3(\xi_3,\tau_3) \tilde{\overline{u_4}}(\xi_4,\tau_4)| \\
&\lesssim \prod_{i=1}^4 \| u_i \|_{X^{0,1/2+}}.
\end{align*}
By Plancherel and H\"older's inequality the left-hand side is bounded by
$$ \| {\cal F}(|\tau_1 - \xi_1^2|^{1/2} \tilde u_1) \|_{L^2_t L^2_x}
\| u_2 \|_{L^6_t L^6_x}
\| u_3 \|_{L^6_t L^6_x}\| u_4 \|_{L^6_t L^6_x}$$
which is acceptable by \eqref{strichartz-6} and the definition of $X^{0,1/2+}$.  This completes the proof of Lemma \ref{a-c}, and hence of Theorem \ref{poly-growth}.

\section{Proof of Theorem \ref{main-2}: first attempt}\label{decomp-sec}

We now prove Theorem \ref{main-2}.  We shall begin by giving a simplified version of the argument which does not capture the full power of $\dist_{H^s}(u_0, \Sigma)$ in \eqref{stab-s}, but introduces the main ideas and establishes the key multi-linear estimates (Lemma \ref{main-est}).  In the next section, we shall give a better version of this argument which recovers most of this power in the next section, and then in the section after that we shall give the most refined version of the argument that gives \eqref{stab-s} with no loss.

This argument is very similar to the previous one, and we have deliberately given the two arguments a nearly identical structure.

\medskip

\divider{Step 0.  Preliminaries; introduction of the modified energy.}

By the existing global well-posedness theory and a standard limiting argument we may assume that $u$ is a global smooth solution which is rapidly decreasing in space.

Fix $0 < s < 1$, $u$, and let $N \gg 1$ be chosen later.
Let $m(\xi)$ be a smooth even real-valued symbol such that $m(\xi) = 1$ for $|\xi| \leq N$ and $m(\xi) = |\xi|^{s-1}/N^{s-1}$ for $|\xi| > 2N$, and let $I$ be the Fourier multiplier
$$ \widehat{If}(\xi) := m(\xi) \hat f(\xi).$$
Thus $I$ is the identity for frequencies $|\xi| \ll N$ and is smoothing of order $1-s$ for high frequencies $|\xi| \gtrsim N$.

In analogy with the quantity \eqref{en-def} used in the proof of Theorem \ref{poly-growth}, we define the \emph{modified energy}\footnote{In a more standard application of the $I$-method, e.g. \cite{ckstt:1}, \cite{ckstt:3}, \cite{ckstt:4}, \cite{keel:mkg}) we would take $E_N(t) = H(Iu(t))$.  The main difference here is thus to add an $L^2$ component to the modified energy in order to make the ground state cylinder an approximate minimizer of the energy.} $E_N(t)$ by
\be{en-def-2}
E_N(t) := L(Iu(t)),
\end{equation}
where the Lyapunov functional $L()$ was defined in \eqref{l-def}.

We now estimate $E_N(0)$.  Write $\sigma := \dist_{H^s}(u_0,\Sigma)$, thus by hypothesis $0 < \sigma \ll 1$ and there exists a ground state $\tilde Q \in \Sigma$ such that
$$ \| u_0 - \tilde Q \|_{H^s} \lesssim \sigma.$$
Applying $I$, we see that
$$ \| Iu_0 - I\tilde Q \|_{H^1} \lesssim N^{1-s} \sigma.$$
On the other hand, since $\tilde Q$ is smooth, its Fourier transform is rapidly decreasing, and
\be{n-0}
\| I \tilde Q - \tilde Q \|_{H^1} \lesssim N^{-C}
\end{equation}
for any $C$.  Thus, if we assume
\be{n-1}
N \gtrsim \sigma^{0-},
\end{equation}
we then have
$$ \| Iu_0 - \tilde Q \|_{H^1} \lesssim N^{1-s} \sigma.$$
By \eqref{weinstein} we thus have
$$ E_N(0) - L(Q) \lesssim N^{2-2s} \sigma^2.$$
We shall make the assumption that
\be{n-s}
N^{2-2s} \sigma^2 \ll 1.
\end{equation}

The heart of the argument shall lie in the following almost conservation law for the modified energy $E_N$.

\begin{lemma}\label{a-c-2}  If $t_0 \in \R$ is such that $E_N(t_0) - L(Q) \ll 1$, then we have
\be{ftoc-2}
E_N(t_0+\delta) = E_N(t_0) + O(N^{-1+})
\end{equation}
where $\delta > 0$ is an absolute constant depending only on $s$.
\end{lemma}

As in the previous section, we do not know if the error bound $O(N^{-1+})$ is sharp; any improvement in this error estimate will ultimately lead to an improvement of the $t^{1-s+}$ factor in \eqref{stab-s} (if one uses the most refined version of the argument below, see Section \ref{final-sec}).

\medskip

\divider{Step 1.  Deduction of a weak form of \eqref{stab-s} from Lemma \ref{a-c-2}.}

If we assume Lemma \ref{a-c-2}, then we may iterate it to obtain
$E_N(t) - L(Q) \lesssim N^{2-2s} \sigma^2$ for all
\be{t-bound}
1 \leq t \ll N^{1-} N^{2-2s} \sigma^2,
\end{equation}

where we assume
\be{n-s-2}
N^{1-} N^{2-2s} \sigma^2 \gg 1
\end{equation}
(note that this automatically implies \eqref{n-1}).

Fix $t$ as above.  Then by \eqref{weinstein} we have\footnote{Strictly speaking, to apply \eqref{weinstein} we must assume beforehand that $Iu(t)$ is close to a ground state.  But by hypothesis this is true at time $t=0$, and so one may proceed by a standard continuity argument which we omit.  This continuity argument will also be used later in the paper without further mention.}
$$ \| Iu(t) - \tilde Q(t) \|_{H^1} \lesssim N^{1-s} \sigma$$
for some ground state $\tilde Q(t) \in \Sigma$ depending on $t$. Using \eqref{n-1} as before, we may modify this to
$$ \| Iu(t) - I\tilde Q(t) \|_{H^1} \lesssim N^{1-s} \sigma$$
which implies that
$$ \| u(t) - \tilde Q(t) \|_{H^s} \lesssim N^{1-s} \sigma$$
and hence that
\be{final-dist}
\dist_{H^s}(u(t), \Sigma) \lesssim N^{1-s} \sigma.
\end{equation}
Optimizing in $N$ subject to \eqref{n-s}, \eqref{t-bound}, \eqref{n-s-2}
we obtain
\be{final-optimum}
\dist_{H^s}(u(t), \Sigma) \lesssim t^{\frac{1-s}{3-2s-}} \sigma^{\frac{1}{3-2s-}}
\end{equation}
whenever $t \ll \sigma^{-1/(1-s)-}$, which is a weak form of \eqref{stab-s}.  In the next two sections we shall use more refined arguments to improve this bound.

It remains to prove Lemma \ref{a-c-2}.  This shall be done in several stages.

\medskip

\divider{Step 2.  Control $u$ at time $t_0$.}

>From the hypothesis $E_N(t_0) - L(Q) \ll 1$ and \eqref{weinstein} (or the Gagliardo-Nirenberg inequality) we have
\be{i-h1}
\| I u(t_0) \|_{H^1} \lesssim 1.
\end{equation}
Up to powers of $N$, this is basically an $H^s$ bound on $u(t_0)$, but to obtain good exponents it is important that we work with the $\| Iu\|_{H^1}$ norm rather than the $\| u \|_{H^s}$.

\medskip

\divider{Step 3.  Control $u$ on the time interval $[t_0-\delta, t_0 + \delta]$.}

We now claim that \eqref{i-h1} implies the spacetime estimate
\be{spacetime-2}
\| Iu \|_{X^{1,1/2+\eps}_{[t_0-\delta,t_0+\delta]}} \lesssim 1
\end{equation}
if $0 < \delta \ll 1$ is a sufficiently small constant and $0 < \eps \ll 1$.  As before, it suffices to show the estimate
\be{junk}
\| I(u_1 \overline u_2 u_3) \|_{1,-1/2+\eps} \lesssim \| Iu_1 \|_{1, 1/2+} \| Iu_2 \|_{1,1/2+} \| Iu_3 \|_{1,1/2+}.
\end{equation}
We can write the left-hand side as
$$ \| I \langle \nabla \rangle(u_1 \overline u_2 u_3) \|_{0,-1/2+\eps},$$
where $I \langle \nabla \rangle$ is the multiplier with symbol $m(\xi) \langle \xi \rangle$.  Since $m(\xi + \eta) \langle \xi+\eta \rangle \lesssim m(\xi) \langle \xi \rangle + m(\eta) \langle \eta \rangle$, we can use the fractional Leibnitz rule, together with \eqref{strichartz-4} and its dual, to obtain this estimate.

\medskip

\divider{Step 4.  Control the increment of the modified energy.}

Now that we have obtained \eqref{spacetime-2}, the next step is again the fundamental theorem of Calculus \eqref{the-ftoc}.  We introduce the nonlinear functional $\Omega(v)$, defined on smooth ($H^2$ will do) functions $v$ on $\R \times \R$ by the formula
\be{omega-def}
\Omega(v) :=
\int_{t_0}^{t_0+\delta}
\left\langle iI\bigl(v(t)_{xx} + F(v(t)) \bigr),
-Iv(t)_{xx} + Iv(t) - F(Iv(t))\right\rangle\ dt.
\end{equation}
>From \eqref{the-ftoc}, \eqref{L-deriv} and \eqref{nls} we have
\be{en-inc}
E_N(t_0+\delta) - E_N(t_0) = 2\Omega(u),
\end{equation}
so in view of \eqref{spacetime-2} it suffices to prove the following estimate, which we shall re-use in later sections.

\begin{lemma}\label{main-est}
We have
$$ |\Omega(v)| \lesssim N^{-1+}$$
whenever
$$\| Iv \|_{X^{1,1/2+\eps}_{[t_0-\delta,t_0+\delta]}} \lesssim 1.$$
\end{lemma}

We now prove Lemma \ref{main-est}.  By repeated integration by parts and a symmetrization (cf. the computations after \eqref{the-ftoc}) we can expand the integrand in \eqref{omega-def} as
\bas
&\langle iIv_{xx}, -F(Iv) \rangle
+ \langle iIF(v), -Iv_{xx} + Iv - F(Iv)\rangle \\
=& \langle iIv_{xx}, -F(Iv) + IF(v) \rangle
+ \langle i F(v), I^2 v \rangle
+ \langle i IF(v), - F(Iv) \rangle \\
=& \Re \bigl( i \Lambda_4(-\xi_1^2 m(\xi_1) (-m(\xi_2) m(\xi_3) m(\xi_4) + m(\xi_2+\xi_3+\xi_4));v,\overline{v},v,\overline{v})\\
&\quad + i \Lambda_4(m(\xi_4)^2,v,\overline{v},v,\overline{v})\\
&\quad -i \Lambda_6(m(\xi_1 + \xi_2 + \xi_3) m(\xi_4) m(\xi_5) m(\xi_6);
v,\overline{v}, v,\overline{v}, v,\overline{v})\bigr)\\
=& \Re \bigl( C_1 \Lambda_4(M'_4; v,\overline{v}, v,\overline{v}) + C_2 \Lambda_4(M''_4; v,\overline{v}, v,\overline{v}) + C_3 \Lambda_6(M_6; v,\overline{v}, v,\overline{v}, v,\overline{v}) \bigr)
\end{align*}
where $C_1$, $C_2$, $C_3$ are imaginary constants and
\bas
 M'_4 &:= \xi_1^2 m(\xi_1) \bigl(m(\xi_2) m(\xi_3) m(\xi_4) - m(\xi_2+\xi_3+\xi_4)\bigr)\\
M''_4 &:= m(\xi_1)^2 - m(\xi_2)^2 + m(\xi_3)^2 - m(\xi_4)^2\\
M_6 &:= m(\xi_1 + \xi_2 + \xi_3) m(\xi_4) m(\xi_5) m(\xi_6) - m(\xi_1) m(\xi_2) m(\xi_3) m(\xi_4 +\xi_5 + \xi_6).
\end{align*}
It thus suffices to show the bounds

\begin{align}
|\int_{t_0}^{t_0+\delta}
\Lambda_4(M'_4; u_1(t), \overline u_2(t), u_3(t), \overline u_4(t))\ dt|
&\lesssim N^{-1} \prod_{i=1}^4 \| Iu_i\|_{1,1/2+} \label{m4}\\
|\int_{t_0}^{t_0+\delta}
\Lambda_4(M''_4; u_1(t), \overline u_2(t), u_3(t), \overline u_4(t))\ dt|
&\lesssim N^{-1} \prod_{i=1}^4 \| Iu_i\|_{1,1/2+} \label{m4'}\\
|\int_{t_0}^{t_0+\delta}
\Lambda_6(M_6; u_1(t), \overline u_2(t), u_3(t), \overline u_4(t), u_5(t), \overline u_6(t))\ dt|
&\lesssim N^{-1} \prod_{i=1}^6 \| Iu_i\|_{1,1/2+}. \label{m6}
\end{align}

As before, we first restrict the Fourier support of $u_i$ to the region $\langle \xi_i \rangle \sim N_i$ for some $N_i \geq 1$, and promise to sum in the $N_i$ later.

 \medskip

\divider{Step 4(a): Proof of \eqref{m4}.}

We first prove \eqref{m4}.  We shall assume that $N_2 \geq N_3 \geq N_4$; the other cases are similar.  We may then assume $N_2 \gtrsim N$ since the symbol vanishes otherwise.

We split into two cases: $N_2 \gg N_3$ and $N_2 \sim N_3$.

\medskip

\divider{Case 4(a).1: $N_2 \gg N_3$.}

We may assume that $N_1 \sim N_2$ since $\xi_1 + \xi_2 + \xi_3 + \xi_4 = 0$.  We then bound the symbol $M'_4$ by $N_1^2 m(N_1) m(N_2) \sim N_1 m(N_1) N_2 m(N_2)$, and estimate this contribution by
$$ N_1 m(N_1) N_2 m(N_2) \int\int |u_1| |u_2| |u_3| |u_4|.$$
Applying Cauchy-Schwartz followed by two applications of \eqref{improved-strichartz} we obtain a bound of
$$ N_1^{1/2} m(N_1) N_2^{1/2} m(N_2) \prod_{i=1}^4 \| u_i \|_{0,1/2+}
\lesssim N_1^{-1/2} N_2^{-1/2} \prod_{i=1}^4 \| Iu_i \|_{1,1/2+}.$$
Summing in the $N_i$ we see that this case is acceptable (we lose some logarithms of $N_1$ and $N_2$ from the $N_3$ and $N_4$ summation, which is why we only end up with a bound of $N^{-1+}$ instead of $N^{-1}$).

\medskip

\divider{Case 4(a).2: $N_2 \sim N_3$.}

Since $\langle \xi_2 + \xi_3 + \xi_4\rangle \sim N_1$, we may bound the symbol $M'_4$ by $N_1^2 m(N_1)^2$, and estimate the contribution by
$$ N_1^2 m(N_1)^2 \int\int |u_1| |u_2| |u_3| |u_4|.$$
We bound $N_1^2 m(N_1)^2 \lesssim N_1 m(N_1) N_2 m(N_2)$ and apply \eqref{strichartz-4} four times to bound this by
$$ N_1 m(N_1) N_2 m(N_2) \prod_{i=1}^4 \| u_i \|_{0,1/2+}
\lesssim \frac{1}{N_3 m(N_3)} \prod_{i=1}^4 \| Iu_i \|_{1,1/2+}.$$
Summing in the $N_i$ we see that this case is acceptable (again, the $N_1$ and $N_4$ summations cost us some logarithms).

\medskip
\divider{Step 4(b): Proof of \eqref{m4'} and \eqref{m6}.}

Now we prove \eqref{m4'}, which is rather easy due to the lack of derivatives in the symbol.  We may assume that $N_1 \geq N_2, N_3, N_4$.  We may assume that $N_1 \gtrsim N$ since the symbol vanishes otherwise. Then if we bound the symbol $M''_4$ by $O(1)$ and use \eqref{strichartz-4} we estimate this term by
$$ \prod_{i=1}^4 \| u_i \|_{0,1/2+} \lesssim \frac{1}{N_1 m(N_1)}
\prod_{i=1}^4 \| Iu_i \|_{1,1/2+}.$$
Summing in the $N_i$ we see that this case is acceptable, again losing some logarithms in the $N_2$, $N_3$, $N_4$ summations.  (Indeed, it is clear one could extract far more decay from this term if desired).  Finally, the estimate \eqref{m6} is similar to \eqref{m4'} (just use \eqref{strichartz-6} instead of \eqref{strichartz-4}).

This concludes the proof of Lemma \ref{main-est}, hence of Lemma \ref{a-c-2}, which gives a weak version of Theorem \ref{main-2}.

\section{Proof of Theorem \ref{main-2}: second attempt}\label{improv-sec}

In the previous section we gave a partial proof of Theorem \ref{main-2}, but with the wrong power of $\dist_{H^s}(u_0, \Sigma)$.  The problem was that we were not really exploiting the fact that $u(t)$ was close to the ground state cylinder $\Sigma$; for instance, the bound \eqref{i-h1} would still be true if $u$ was at a distance $\sim 1$ from the ground state cylinder in the $\| I u \|_{H^1}$ metric.  To improve upon these results we must consider not just $u(t)$, but also the difference $w(t) := u(t) - Q(t)$ between $u$ and an appropriate ground state $Q(t) \in \Sigma$.  In particular we wish to exploit the fact that $w$ has small norm (with a bound which depends linearly on $\sigma$).

In this section we use the above ideas to refine the argument of the previous section, and obtain a near miss to \eqref{stab-s}.  Unfortunately our power of $\dist_{H^s}(u_0, \Sigma)$ will still be off by an epsilon, mainly because the smoothing operator $I$ does not quite preserve the ground state cylinder (so that $IQ(t)$ is not a ground state).  To fix this problem and get the sharp power of $\dist_{H^s}(u_0, \Sigma)$ requires some further refinements which we delay until the next section in order to simplify the exposition.

We now turn to our second attempt at proving Theorem \ref{main-2}, deliberately repeating much of the structure of the arguments from previous sections.

\medskip

\divider{Step 0.  Preliminaries; introduction of the modified energy.}

We make the same reductions as the previous section, and leave the definition of the modified energy $E_N(t)$ from \eqref{en-def-2} unchanged.  The main difference is that we sharpen Lemma \ref{a-c-2} to

\begin{lemma}\label{a-c-3}  If $t_0 \in \R$ is such that $E_N(t_0) \leq L(Q) + \tilde \sigma^2$ for some $N^{-C} < \tilde \sigma \ll 1$ for some arbitrary constant $C$, then we have
\be{ftoc-3}
E_N(t_0+\delta) = E_N(t_0) + O(N^{-1+} \tilde \sigma^2)
\end{equation}
where $\delta > 0$ is an absolute constant depending only on $s$.
\end{lemma}

As in the previous section, we do not know if the factor $O(N^{-1+})$ can be improved.  The quadratic exponent $\tilde \sigma^2$ is probably sharp, however, since the derivative of $E_N(t_0)$ will contain terms which are quadratic in the difference $w(t) = u(t) - Q(t)$.

\medskip

\divider{Step 1.  Deduction of a slightly weakened form of \eqref{stab-s} from Lemma \ref{a-c-3}.}

Assume Lemma \ref{a-c-3} for the moment.  We substitute this lemma, with $\tilde \sigma \sim N^{1-s} \sigma$, in place of Lemma \ref{a-c-2} in Step 1 of the previous section; note that $\tilde \sigma$ is admissible by \eqref{n-1}.  The argument then proceeds as before, with \eqref{t-bound} replaced by
$$
1 \leq t \ll N
$$
(and \eqref{n-s-2} is no longer needed).  For this range of $t$ we obtain \eqref{final-dist}.  Optimizing in $N$ subject to \eqref{n-1}, \eqref{n-s}, \eqref{t-bound} we find that we may improve \eqref{final-optimum} to
\be{final-optimum-2}
\dist_{H^s}(u(t), \Sigma) \lesssim t^{1-s+} \sigma^{1-}
\end{equation}
whenever
$$ t \ll \sigma^{-\frac{1}{1-s}-}.$$
This is within an epsilon of \eqref{stab-s}.  To remove this last epsilon we shall need a more refined argument, presented in the next section.

It remains to prove Lemma \ref{a-c-3}.  This shall be done in the usual sequence of stages.  The main difference is that, instead of controlling $u$, we shall control a difference $w(t) = u(t) - Q(t)$ between $u(t)$ and a suitably chosen ground state $Q(t)$.  This will let us recover the powers of $\tilde \sigma$ in \eqref{ftoc-3}.

\medskip

\divider{Step 2.  Control $w$ at time $t_0$.}

Fix $t_0$.  By \eqref{weinstein} there exists a ground state $Q_{t_0} \in \Sigma$ depending on $t_0$ such that
$$ \| Iu(t_0) - Q_{t_0} \|_{H^1} \lesssim \tilde \sigma.$$
The above estimate asserts that $u(t_0)$ is in some sense close to $Q_{t_0}$.  If $u(t_0)$ was in fact equal to $Q_{t_0}$, then the evolution of $u(t)$ would follow the curve of ground states\footnote{This heuristic is a little inaccurate because of the presence of the $I$, and also because the error $Iu(t_0) - Q_{t_0}$ can affect the modulation parameters of the approximating ground state.  We address these issues in the next section, when we remove the epsilon losses which arise from the arguments in this section.} $Q(t): \R \to \Sigma$ defined by
\be{curve}
Q(t) := e^{i(t-t_0)} Q_{t_0}.
\end{equation}
Thus it is natural to define $w(t) := u(t) - Q(t)$.  By \eqref{n-0} and the previous we have
\be{w-init}
\| I w(t_0) \|_{H^1} \lesssim \tilde\sigma,
\end{equation}
where we have used the hypothesis that $\tilde \sigma \gtrsim N^{-C}$ for some $C$.

\medskip

\divider{Step 3.  Control $w$ on the time interval $[t_0-\delta, t_0 + \delta]$.}

The idea is now to run a local well posedness argument for $w$ instead of $u$ in order to gain the powers of $\tilde \sigma$ in \eqref{ftoc-3}.  Specifically, we wish to obtain the spacetime estimate
\be{spacetime-3}
\| Iw \|_{X^{1,1/2+\eps}_{[t_0-\delta,t_0+\delta]}} \lesssim \tilde \sigma
\end{equation}
for some small absolute constants $0 < \delta \ll 1$ and $0 < \eps \ll 1$.

To obtain this, we observe from \eqref{nls} and \eqref{Q-def} that $w$ obeys the difference equation
\be{w-diff}
iw_t + w_{xx} = -G(w(t),Q(t))
\end{equation}
where $G$ is the nonlinear expression
$$
G(w(t),Q(t)) := F(Q(t)+w(t)) - F(Q(t)).$$
The exact form of $G$ is not important, save for the fact that $G$ is cubic in $Q(t)$, $\overline{Q(t)}$, $w$, $\overline{w}$, and is always at least linear in $w$, $\overline{w}$.
By \eqref{w-init}, \eqref{w-diff} and the standard $X^{s,b}$ energy estimate (as in previous sections), we thus have
$$ \| Iw \|_{X^{1,1/2+\eps}_{[t_0-\delta,t_0+\delta]}} \lesssim \tilde \sigma
+ \delta^{0+} \| IG(w) \|_{X^{1,-1/2+2\eps}_{[t_0-\delta,t_0+\delta]}}.$$
By several applications of \eqref{junk}, and the observation that $Q(t)$, $\overline{Q(t)}$ are Schwartz in $x$ and thus live in every (time-localized) $X^{s,b}$ space, we therefore have
$$ \| Iw \|_{X^{1,1/2+\eps}_{[t_0-\delta,t_0+\delta]}} \lesssim \tilde \sigma
+ \delta^{0+} \| Iw\|_{X^{1,1/2+\eps}_{[t_0-\delta,t_0+\delta]}} + \delta^{0+} \| Iw\|_{X^{1,1/2+\eps}_{[t_0-\delta,t_0+\delta]}}^3.$$
The claim \eqref{spacetime-3} then follows if $\delta$ is small enough by standard continuity (or iteration) arguments.

\medskip

\divider{Step 4.  Control the increment of the modified energy.}

We now prove \eqref{ftoc-3}.  By \eqref{en-inc} it suffices to show that
$$ \Omega(Q(t) + w(t)) = O(N^{-1+} \tilde \sigma^2).$$
We could do this by direct computation, but we present instead a simple argument (based on isolating the terms in $\Omega(Q(t) + w(t))$ which are linear in $w$) which obtains this bound as a nearly automatic consequence of Lemma \ref{main-est}.

To prove the above estimate it suffices by the hypothesis $\tilde \sigma \gtrsim N^{-C}$ to prove the more general bound
\be{k-grow}
\Omega(Q(t) + k \frac{w(t)}{\tilde \sigma}) = O(N^{-1+} |k|^2) + O(N^{-C-1} |k|)
\end{equation}
for any real number $k$ such that $|k| \lesssim 1$.

To prove \eqref{k-grow}, first observe from \eqref{omega-def} that the left-hand side of \eqref{k-grow} is a polynomial $P(k)$ in $k$ of degree at most 6 (with the coefficients depending on $Q(t)$, $w(t)$, $\sigma$, of course).  From \eqref{omega-def}, \eqref{Q-def}, and integration by parts we see that the constant term of $P$ is zero.  Also, from Lemma \ref{main-est} and \eqref{spacetime-3} we see that the left-hand side of \eqref{k-grow} is $O(N^{-1+})$ for all $|k| \lesssim 1$.  Thus all the coefficients of $P(k)$ are $O(N^{-1+})$.  To finish the argument it suffices to show that the linear term of $P(k)$ is $O(N^{-C-1})$.  Equivalently, it suffices to show that the portion of $\Omega(Q(t)+w(t))$ which is linear in $w$, $\overline{w}$ is $O(N^{-C-1} \tilde \sigma)$.

To show this we return to \eqref{en-inc}, and exploit the fact that $IQ(t)$ is nearly a minimizer of $L$.  If we compute $E_N(t) = L(IQ(t) + Iw(t))$ using \eqref{l-def} and suppress terms which are not linear in $w$, $\overline{w}$, we obtain

\bas
E_N(t) = &
2 \langle IQ(t)_x, Iw(t)_x \rangle + 2 \langle IQ(t), Iw(t)\rangle \\
& \quad - 2
\langle F(IQ(t)), Iw(t) \rangle + \hbox{nonlinear terms}
\end{align*}
Integrating by parts and using \eqref{Q-def}, we obtain
$$
E_N(t) =
2 \langle w(t), R(t) \rangle + \hbox{nonlinear terms}$$
where
$$ R(t) := I(IF(Q(t)) - F(IQ(t))).$$
Differentiating this in time using \eqref{w-diff}, and again suppressing nonlinear terms, we obtain
\bas
\partial_t E_N(t) =
&2 \bigl\langle i w_{xx}(t) + 2i w(t) \overline{Q(t)} Q(t) + i Q(t) \overline{w}(t)Q(t), R(t) \bigr\rangle \\
&\quad +
2 \langle w(t), \partial_t R(t) \rangle +
\hbox{nonlinear terms}
\end{align*}
We can bound the linear terms by (for instance)
$$ \lesssim \| w(t) \|_2 (\| R(t) \|_{H^2} + \| \partial_t R(t) \|_2).$$
Because $Q(t)$ is Schwartz, the operator $I$ is almost the identity on $Q(t)$ or $F(Q(t))$, and so it is easy to see that\footnote{Note that the time derivative $\partial_t$ on $R(t)$ is not dangerous because the time-dependence of $Q(t)$ and $R(t)$ are given by a simple phase rotation $e^{it}$.  In the next section however we shall need a more sophisticated time evolution for the approximating ground state $Q(t)$.}
$$ \| R(t) \|_{H^2} + \| \partial_t R(t) \|_2 \lesssim N^{-C-1}.$$
Also from \eqref{spacetime-3} we have $\| w(t)\|_2 \lesssim \tilde \sigma$.  Thus the linear part of $\Omega(Q(t)+w(t))$ is $O(N^{-C-1} \tilde \sigma)$ as desired.
This proves \eqref{ftoc-3}, which then almost gives \eqref{stab-s}.

\section{Proof of Theorem \ref{main-2}: final argument}\label{final-sec}

In the previous sections one had to assume \eqref{n-1} in order to use the approximation $Q \approx IQ$.  Ultimately this assumption caused us to miss \eqref{stab-s} by an epsilon.  In order to avoid this loss we shall need to avoid using the approximation $Q \approx IQ$, at least in situations in which one does not have enough powers of $\sigma$ in the estimates.

\medskip

\divider{Step 0.  Preliminaries; introduction of the modified energy.}

This step is the same as in the previous section, except that
we need to replace the modified energy $L(Iu)$ from \eqref{en-def-2} by a slightly different quantity (because $IQ$ is not quite a minimizer of $L$).

To motivate the argument we first recall the more precise statement
of Weinstein's estimate \eqref{weinstein}:

\begin{lemma}\label{wlemma}\cite{wein}, \cite{wein:modulate}  Let
$\tilde Q \in \Sigma$ be a ground state, and let $w \in H^1$ obey the
two orthogonality conditions\footnote{If one writes $u := \tilde Q +
w$, then these orthogonality conditions have the geometric
interpretation that $\tilde Q$ is the closest ground state in
$\Sigma$ to $u$, as measured in $H^1$ norm.}
$$ \langle w, A F(\tilde Q) \rangle = 0 \hbox{ for } A = i,
\partial_x$$
Then, if the $H^1$ norm of $w$ is sufficiently small, we have the
coercivity estimate
$$ L(\tilde Q + w) - L(\tilde Q) = L(\tilde Q + w) - L(Q) \sim \| w
\|_{H^1}^2.$$
\end{lemma}

Note that the two anti-selfadjoint operators $A = i, \partial_x$ are
the infinitesimal generators of the phase rotation and translation
groups respectively, both of which preserve the ground state cylinder
$\Sigma$.

We shall need to apply this lemma to estimate quantities of the form
$L(Q(t)+Iw(t))$ (which will be our substitute for $L(Iu)$).  Thus we
shall require $w$ to obey the orthogonality conditions\footnote{This
trick of choosing $Q(t)$ to obey carefully selected orthogonality
conditions is very common in stability analysis, see e.g.
\cite{merle}.}
\be{ortho}
\langle w(t), A IF(Q(t)) \rangle = 0 \hbox{ for } A = i, \partial_x
\end{equation}
for all times $t$.

To use Lemma \ref{wlemma} we thus need a decomposition $u = Q(t) +
w(t)$ of the solution $u$ which obeys \eqref{ortho}.  This is the
purpose of the following Lemma.

\begin{lemma}\label{approx}  If $u \in H^s(\R)$ is such that
$\dist_{H^s}(u, \Sigma) \ll N^{s-1}$, and $N$ is sufficiently large
depending on $s$, then we may decompose $u = \tilde Q + w$, where
$\tilde Q \in \Sigma$ is a ground state,
$w$ obeys the orthogonality conditions
\be{ortho-tilde}
\langle w, AIF(\tilde Q) \rangle = 0 \hbox{ for } A = i, \partial_x
\end{equation}
and the bound
\be{shoop-alt}
\| Iw\|_{H^1} \lesssim N^{1-s} \dist_{H^s}(u,\Sigma).
\end{equation}
\end{lemma}

\begin{proof}  Define the metric $d$ on $H^s$ by
$$ d(u,v) := \| I(u-v)\|_{H^1} = \left(\langle I(u-v), I(u-v) \rangle
+ \langle \partial_x I(u-v), \partial_x I(u-v) \rangle\right)^{1/2}.$$
Clearly we have
$$ d(u,\Sigma) \lesssim N^{1-s} \dist_{H^s}(u,\Sigma) \ll 1.$$
We now claim that there exists a ground state $Q'$ in $\Sigma$ which
minimizes $d(u,Q')$, so that $d(u,Q') = d(u,\Sigma)$.  To see this,
first observe (since $u$ was assumed to be smooth and rapidly
decreasing) that $\langle Iu, I e^{i\theta} Q(x-x_0) \rangle \to 0$
as $x_0 \to \pm \infty$, and similarly for
$\langle \partial_x Iu, \partial_x I e^{i\theta} Q(x-x_0) \rangle$.
By orthogonality we thus have
$$ d(u, e^{i\theta} Q(x-x_0)) \to (\| Iu \|_{H^1}^2 + \|
IQ\|_{H^1}^2)^{1/2} \gtrsim 1 \gg d(u, \Sigma)$$
as $x_0 \to \pm \infty$.  Thus, in order to minimize $d(u, Q')$ on
$\Sigma$, it suffices to restrict $Q'$ to a compact subset of the
cylinder $\Sigma$.  By compactness and smoothness we then see that a
minimizer $Q'$ must exist.

Observe that the statement and conclusion of Lemma \ref{approx} is
invariant under translations and modulations of $u$ (and hence of
$\tilde Q$ and $w$).  By these invariances we may thus assume that
the minimizer $Q'$ is attained at $Q' = Q$, thus
\be{shoop}
d(u,Q) = d(u, \Sigma) \lesssim N^{1-s} \dist_{H^s}(u,\Sigma) \ll 1.
\end{equation}

The tangent space of $\Sigma$ at $Q$ is spanned by $iQ$ and $Q_x$.
Differentiating
$$d(u,Q)^2 = \langle I(u-Q), I(u-Q) \rangle + \langle \partial_x
I(u-Q), \partial_x I(u-Q) \rangle$$
in these directions, we thus see that
\bas \langle I(u-Q), iIQ \rangle + \langle \partial_x I(u-Q), iI Q
\rangle &= 0; \\
\langle I(u-Q), IQ_x \rangle + \langle \partial_x I(u-Q), \partial_x
IQ_x
\rangle &= 0.
\end{align*}
Integrating by parts and using \eqref{Q-def} we obtain
\be{near-ortho}
\langle\tilde w, A I^2 F(Q) \rangle = 0 \hbox{ for } A = i, \partial_x
\end{equation}
where $\tilde w := u-Q$.  This is almost what we want, except that we have
$I^2$ instead of $I$.  To rectify this we shall need to use
perturbation theory to shift and modulate the ground state slightly.
In other words, we set $\tilde Q := e^{i\theta} Q(x-x_0)$ for some
$|\theta|, |x_0| \ll 1$ to be chosen later.  Writing $q := \tilde Q -
Q$ and
$$w := u - \tilde Q = \tilde w - q,$$
we see from \eqref{near-ortho} that \eqref{ortho-tilde} becomes
$$ \langle \tilde w - q, A I F(Q+q) \rangle - \langle \tilde w, A I^2 F(Q)
\rangle = 0 \hbox{ for } A = i, \partial_x.$$
We rearrange this as
\be{system} \langle q, AIF(Q+q) \rangle - \langle \tilde w,
AI(F(Q+q)-F(Q)) \rangle =
\langle I\tilde w, A(1-I) F(Q) \rangle \hbox{ for } A = i,
\partial_x.
\end{equation}
Since $F(Q)$ is Schwartz, the right-hand side is $O(N^{-100}
\|I\tilde w\|_{H^1}) = O(N^{-99} \dist_{H^s}(u,\Sigma))$
by \eqref{shoop}.  Now
we expand the left-hand side to first order in $\theta$, $\xi_0$.
Observe that
$$ q = \tilde Q-Q = \theta iQ - x_0 Q_x + O_{H^2}(|\theta|^2 +
|x_0|^2)$$
where $O_{H^2}(X)$ denotes a quantity whose $H^2$ norm is $O(X)$.
Similarly we have
$$ F(Q+q)-F(Q) = \theta G + x_0 H + O_{H^2}(|\theta|^2 + |x_0|^2)$$
where $G$, $H$ are explicit Schwartz functions whose exact form is
not important for us.  Thus \eqref{system} becomes
\begin{align*}
 \theta \langle iQ, AIF(Q) \rangle &- x_0 \langle Q_x, AIF(Q) \rangle
- \theta \langle \tilde w, AIG \rangle - x_0 \langle \tilde w, AIH \rangle\\
&= O(N^{-99} \dist_H^s(u, \Sigma)) + O(|\theta|^2 + |x_0|^2) \hbox{
for } A = i, \partial_x,
\end{align*}
which we can write as a matrix system
\be{system-2}
\begin{split}
&\left( \begin{array}{ll}
\langle iQ, iIF(Q) \rangle - \langle \tilde w, iIG \rangle &
-\langle Q_x iIF(Q) \rangle - \langle \tilde w, iIH \rangle \\
\langle iQ, \partial_x IF(Q) \rangle - \langle \tilde w, \partial_x IG
\rangle &
-\langle Q_x \partial_x IF(Q) \rangle - \langle \tilde w,
\partial_x IH
\rangle \\
\end{array}
\right)
\left( \begin{array}{l}
\theta \\
x_0
\end{array}
\right)\\
&\qquad\qquad = O(N^{-99} \dist_H^s(u, \Sigma)) + O(|\theta|^2 +
|x_0|^2).
\end{split}
\end{equation}
Since $G$ and $H$ are Schwartz, we have
$$ \langle \tilde w, AIG \rangle, \langle \tilde w, AIH \rangle = O(\|
I\tilde w
\|_{H^1})
= O(d(u,Q)) \ll 1 \hbox{ for } A = i, \partial_x$$
by \eqref{shoop}.
Also, an easy integration by parts using \eqref{n-0} gives the
estimates
\bas \langle iQ, iI F(Q) \rangle &= \|Q\|_4^4 + O(N^{-100}) \\
\langle Q_x, iI F(Q) \rangle,
\langle iQ, \partial_x I F(Q) \rangle &= O(N^{-100})\\
\quad \langle Q_x, \partial_x I F(Q) \rangle &= 3 \int Q^2 Q_x^2 +
O(N^{-100}).
\end{align*}
Thus if $N$ is large enough, the matrix in \eqref{system-2} has a
non-degenerate Jacobian at the point $(\theta, x_0) = (0,0)$.  By the
inverse function theorem (or the contraction mapping theorem) we can
thus solve \eqref{system} for some $\theta, x_0 = O(N^{-99}
\dist_{H^s}(u,\Sigma))$.  The condition \eqref{shoop-alt} then
follows from \eqref{shoop} (since the distance between $Q$ and
$\tilde Q$ is $O(|\theta|+|x_0|) = O(N^{-99} \dist_{H^s}(u,\Sigma))$
in any reasonable norm).
\end{proof}

Applying this Lemma at each time $t$ we thus have ground states
$Q(t)$ and a function $w(x,t)$ such that $u(t) = Q(t) + w(t)$ and the
orthogonality relations \eqref{ortho} hold for all times $t$ for
which $\dist_{H^s}(u, \Sigma) \ll N^{s-1}$.  (We can ignore this
latter condition by standard continuity arguments, since we will
eventually verify this hypothesis (with some room to spare) at the
end of the argument).

We now redefine the modified energy $E_N(t)$ as
\be{en-def-3}
E_N(t) := L(Q(t) + Iw(t)).
\end{equation}
This should be compared to the energy $L(IQ(t) + Iw(t))$ from \eqref{en-def-2} used in the previous section.  From Lemma \ref{wlemma} we have the analogue of \eqref{weinstein}
\be{tilde-en}
E_N(t) - L(Q) \sim \| Iw(t)\|_{H^1}^2
\end{equation}
whenever $\| Iw(t)\|_{H^1}$ is sufficiently small.  In particular, at time zero we see from Lemma \ref{approx} that
$$ E_N(0) - L(Q) \sim \| Iw(0)\|_{H^1}^2
\lesssim N^{2-2s} \sigma^2.$$

The analogue of Lemma \ref{a-c-3} is

\begin{lemma}\label{a-c-4}  If $t_0 \in \R$ is such that $E_N(t_0) \leq L(Q) + \tilde \sigma^2$ for some $0 < \tilde \sigma \ll 1$, then \eqref{ftoc-3} holds for some $\delta > 0$ which is an absolute constant depending only on $s$.
\end{lemma}

\medskip

\divider{Step 1.  Deduction of \eqref{stab-s} from Lemma \ref{a-c-4}.}

We repeat Step 1 of the previous section, using \eqref{tilde-en} as a substitute for \eqref{weinstein}.  The main differences are that there is no hypothesis of the form $\tilde \sigma > N^{-C}$ in Lemma \ref{a-c-4}, so that there is no need to assume \eqref{n-1}.  We omit the details.

Note that while there is some logarithmic losses in the $O(N^{-1+})$ factor of \eqref{ftoc-3}, there are none in the $\tilde \sigma^2$ factor.  This is what allows us to get the sharp power of $\dist_{H^s}(u_0, \Sigma)$ in \eqref{stab-s}, although the power $t^{1-s+}$ is probably not sharp.

It remains to prove Lemma \ref{a-c-4}, which we divide into the now-familiar sequence of steps.

\medskip

\divider{Step 2.  Control $w$ at time $t_0$.}

We now prove Lemma \ref{a-c-4}.  Fix $t_0$.  By \eqref{tilde-en} we again obtain \eqref{w-init}.

\medskip

\divider{Step 3.  Control $w$ on the time interval $[t_0-\delta, t_0 + \delta]$.}

As before, the next step is to obtain \eqref{spacetime-3}.
It will be convenient to write $Q(t)$ more explicitly as
\be{qt-def}
Q(x,t) =: e^{i\theta(t)} e^{it} Q(x - x_0(t))
\end{equation}
where the modulation parameters $\theta(t)$, $x_0(t)$ are real-valued; we shall shortly derive equations for the evolution of these parameters.  This ansatz should be compared with \eqref{curve}; note that the standard ground state $u(t) = e^{it} Q$ occurs when $\theta, x_0, w$ all identically vanish.

>From \eqref{nls}, \eqref{Q-def}, \eqref{qt-def} we see that $w$ obeys the difference equation
\be{w-evolve}
\begin{split}
iw_t + w_{xx} &= -G(w(t),Q(t)) -Q(t)-iQ(t)_t \\
&= -G(w(t),Q(t)) + \dot \theta(t) Q(t) + i \dot x_0(t) \partial_x Q(t)
\end{split}
\end{equation}
where $G$ was defined in the previous section.  We may thus repeat Step 3 of the previous section, provided that we can show the estimate
$$ \| \dot \theta(t) Q(t) + i \dot x_0(t) \partial_x Q(t) \|_{X^{1,-1/2+\eps}_{[t_0-\delta,t_0+\delta]}} \lesssim
\| Iw \|_{X^{1,1/2+}_{[t_0-\delta,t_0+\delta]}}.$$
Since $Q(t)$ is Schwartz, it suffices to show that
\be{theta-crew}
|\dot \theta(t)|, |\dot x_0(t)| \lesssim \|w(t)\|_{H^s} \lesssim \| Iw(t) \|_{H^1}.
\end{equation}

To do this, we argue in a manner reminiscent of the proof of Lemma \ref{approx}.  We first introduce the renormalized function $\tilde w(t)$ defined by
$$ \tilde w(x,t) := e^{-i\theta(t)} e^{-it} w(x + x_0(t)).$$
>From \eqref{w-evolve} and \eqref{qt-def} we see that $\tilde w$ evolves according to the equation
\be{w-evolve-2}
i\tilde w_t + \tilde w_{xx} = -G(\tilde w,Q) + \dot \theta(t) (Q+\tilde w)
+ \tilde w + \dot x_0(t) \partial_x (Q+\tilde w).
\end{equation}
On the other hand, from \eqref{ortho} we have the orthogonality relations
$$\langle \tilde w(t), A_j \rangle = 0 \hbox{ for } j = 0, 1,$$
where $A_0$, $A_1$ are the Schwartz functions
$$ A_0 := I(iQ^3); \quad A_1 := I\partial_x(Q^3).$$
Let $j=0,1$.  Differentiating the previous in time and then applying \eqref{w-evolve-2} we obtain
$$
- \langle i\tilde w_{xx}, A_j \rangle = -\dot \theta(t) \langle i(Q+\tilde w), A_j \rangle
+ \langle i\tilde w, A_j \rangle + \dot x_0(t) \langle \partial_x (Q+\tilde w) A_j \rangle - \langle iG(\tilde w, Q), A_j \rangle$$
which we rewrite as
$$ \langle i(Q+\tilde w), A_j \rangle \dot \theta(t)
- \langle \partial_x (Q+\tilde w), A_j \rangle \dot x_0(t)
= -\langle \tilde w, i (A_j)_{xx} \rangle + \langle G(\tilde w,Q), iA_j \rangle$$
This is a linear system of two equations $j=0,1$ and two unknowns $\dot \theta(t)$, $\dot x_0(t)$.  To invert this system we first observe that the right
hand side is $O(\|w\|_{H^s})$ since the $A_j$ are Schwartz.  Also, since $IQ - Q$ has norm $O(N^{-100})$ in any reasonable space, we have the coefficient estimates
\bas
 \langle i(Q+\tilde w), A_0 \rangle &= \langle iQ, iQ^3 \rangle + \langle i (IQ-Q), iQ^3 \rangle + O(\|\tilde w\|_{H^s}) \\ &= \|Q\|_4^4 + O(N^{-100}) + O(\|\tilde w\|_{H^s})\\
\langle i(Q+\tilde w), A_1 \rangle &= \langle i\tilde w, A_1 \rangle + O(N^{-100})\\
&= O(\|\tilde w\|_{H^s}) + O(N^{-100})\\
\langle \partial_x (Q+\tilde w), A_0 \rangle &= -\langle \tilde w, \partial_x A_0 \rangle + O(N^{-100})\\
&= O(\|\tilde w\|_{H^s}) + O(N^{-100})\\
\langle \partial_x (Q+\tilde w), A_1 \rangle &=
\langle \partial_x Q, \partial_x (Q^3) \rangle + \langle \partial_x (IQ-Q),
\partial_x (Q^3) - \langle w, \partial_x A_1 \rangle
\\ &= - \langle Q_{xx}, Q^3 \rangle
+ O(N^{-100}) + O(\|\tilde w\|_{H^s}).
\end{align*}
Observe that the absolute constants $\|Q\|_4^4$ and $-\langle Q_{xx}, Q^3 \rangle = 3 \int Q^2 Q_x^2$ are both non-zero.  Thus if $N$ is sufficiently large and $\|\tilde w\|_{H^s} = \| w \|_{H^s}$ is sufficiently small, we can invert the above $2 \times 2$ linear system and obtain the desired bounds \eqref{theta-crew}.
Note that this argument shows that one can make $\theta(t)$ and $x_0(t)$ differentiable in $t$.

This completes the proof of \eqref{spacetime-3}.  In particular we have the bounds
\be{w-energy}
\| w(t) \|_{H^s} \lesssim \| Iw(t)\|_{H^1} \lesssim \sigma
\end{equation}
on the interval $[t_0-\delta, t_0+\delta]$, so from \eqref{theta-crew} we have
\be{theta-2}
|\dot \theta(t)|, |\dot x_0(t)| \lesssim \tilde \sigma.
\end{equation}

\medskip

\divider{Step 4.  Control the increment of the modified energy.}

We now prove \eqref{ftoc-3} again (or more precisely, we prove \eqref{ftoc-3} for the new definition \eqref{en-def-3} of the modified energy).  We write
$$ E_N(t) = L(e^{-it} Q(t) + e^{-it} Iw(t)).$$
We now use \eqref{L-deriv} to obtain
$$
\partial_t E_N(t)
= 2\langle Q(t)_t - iQ(t) + Iw_t - iIw, -\partial_{xx}(Q(t)+Iw) + (Q(t)+Iw) - F(Q(t)+Iw) \rangle.
$$
We simplify the right factor using \eqref{Q-def} to obtain
$$
\partial_t E_N(t)
= 2\langle Q(t)_t - iQ(t) + Iw_t - iIw, -Iw_{xx} + Iw - G(Iw,Q(t)) \rangle.
$$

>From \eqref{qt-def}, \eqref{w-evolve} we have
\bas Q(t)_t - iQ(t) + Iw_t - iIw =&\ i \dot \theta(t) (Q(t)-IQ(t))- \dot x_0(t) \partial_x (Q(t)-IQ(t)) \\&+ I(iw_{xx} - iw + iG(w,Q(t))).
\end{align*}
To show \eqref{ftoc-3}, it thus suffices by \eqref{the-ftoc}, \eqref{theta-2} to show the bounds
\begin{align}
\int_{t_0}^{t_0+\delta}
\tilde \sigma | \langle Q(t) - IQ(t), -Iw_{xx} + Iw - G(Iw,Q(t)) \rangle |\ dt
& \lesssim N^{-1+} \tilde \sigma^2 \label{a1}\\
\int_{t_0}^{t_0+\delta}
\tilde \sigma | \langle \partial_x(Q(t) - IQ(t)), -Iw_{xx} + Iw - G(Iw,Q(t)) \rangle |\ dt
& \lesssim N^{-1+} \tilde \sigma^2 \label{a2}\\
|\int_{t_0}^{t_0+\delta}
 \langle I(iw_{xx} - iw + iG(w,Q(t))), -Iw_{xx} + Iw - G(Iw,Q(t)) \rangle \ dt|
& \lesssim N^{-1+} \tilde \sigma^2. \label{b}
\end{align}

To prove \eqref{a1}, \eqref{a2} we use integration by parts to move all the derivatives onto $Q(t)-IQ(t)$.  This function has a norm of $N^{-100}$ in any reasonable space.  By \eqref{w-energy} we thus see that these terms are acceptable.

We now prove \eqref{b}.  Again we could do this by direct computation, but we shall instead just use the work that we have already done in previous sections.  By \eqref{spacetime-3} it suffices to show
$$ |\int_{t_0}^{t_0+\delta}
\langle iIw_{xx}- iIw + iIG(w, Q(t)), -Iw_{xx} + Iw - G(Iw,Q(t)) \rangle| \lesssim
N^{-1+} \| Iw \|_{X^{1,1/2+\eps}_{[t_0-\delta,t_0+\delta]}}^2.$$
The left-hand side consists of multilinear expressions of order between 2 and 6 in $w$.  Thus (cf. Step 4 of the previous section) it suffices to prove the estimate
$$ |\int_{t_0}^{t_0+\delta}
\langle iI\tilde w_{xx} - iI\tilde w + iIG(\tilde w, Q(t)), -I\tilde w_{xx} + I\tilde w - G(I\tilde w,Q(t)) \rangle| \lesssim
N^{-1+}$$
for all spacetime functions $\tilde w$ in the unit ball of $X^{1,1/2+\eps}_{[t_0-\delta,t_0+\delta]}.$  Note that the parameter $\tilde \sigma$ has totally disappeared, and so we are now able to lose factors of $O(N^{-100})$ as necessary.

We now convert the above expression into one which can be dealt with by Lemma \ref{main-est}.  Fix $\tilde w$, and define the functions $v(t)$ by
$$ v(t) := \tilde w(t) + Q(t).$$
>From \eqref{Q-def} we have
$$ iI\tilde w_{xx} - iI\tilde w + iIG(\tilde w, Q(t))
= I(i v_{xx} - iv +  i F(v)).$$
Similarly we have
$$
-I\tilde w_{xx} + I\tilde w - G(I\tilde w,Q(t))
= -I v_{xx} + Iv - F(Iv) + B(t)$$
where the error $B(t)$ is given by
$$ B(t) := [F(I\tilde w + IQ(t)) - F(I\tilde w + Q(t))]
- [IF(Q(t)) - F(Q(t))].$$
Recall that $IQ(t)$ is within $O(N^{-100})$ to $Q(t)$ in $H^1$ norm (say), and similarly for $IF(Q(t))$ and $F(Q(t))$.  Since $I\tilde w$ is also bounded in $H^1$, we thus have
$$ \| B(t) \|_{H^1} \lesssim N^{-100},$$
and so the contribution of $B(t)$ is easily seen to be acceptable.  Thus it remains to show that
$$ |\int_{t_0}^{t_0+\delta}
\langle I(iv_{xx} + iv - iF(v), -Iv_{xx} + Iv - F(Iv) \rangle| \lesssim
N^{-1+}.$$
>From integration by parts we have the identity
$$ \langle iIv, -Iv_{xx} + Iv - F(Iv) \rangle = 0$$
so it suffices to show
$$ |\int_{t_0}^{t_0+\delta}
\langle I(iv_{xx} - iF(v), -Iv_{xx} + Iv - F(Iv) \rangle| \lesssim
N^{-1+}.$$
But this is immediate from Lemma \ref{main-est} and \eqref{omega-def}.
This finishes the proof of \eqref{ftoc-3}, and thus Theorem
 \ref{main-2} is (finally!) completely proved.

\medskip

\noindent
 Email Addresses:\\
    {colliand\verb+@+math.toronto.edu};
 {keel\verb+@+math.umn.edu};
   {gigliola\verb+@+math.mit.edu}\\
    {takaoka\verb+@+math.sci.hokudai.ac.jp};
    {tao\verb+@+math.ucla.edu}
 \medskip

   Received May 2002; revised October 2002.

   \medskip

\end{document}